\documentstyle[12pt]{article}

\begin{document}

\author{S.V. Ludkovsky}
\title{Infinitely divisible distributions over locally compact
non-archimedean fields}
\date{25.06.2007}
\maketitle

\begin{abstract}
The article is devoted to stochastic processes with values in
finite-dimensional vector spaces over infinite locally compact
fields with non-trivial non-archimedean valuations. Infinitely
divisible distributions are investigated. Theorems about their
characteristic functionals are proved. Particular cases are
demonstrated.
\end{abstract}

\section{Introduction}
\par It is well-known that infinitely divisible
distributions play very important role in the theory of stochastic
processes over fields of real and complex numbers
\cite{dal,gihsko,petrov}. But over infinite fields with
non-archimedean non-trivial norms they were not studied. This
article is devoted to  infinitely divisible distributions of
stochastic processes in vector spaces over locally compact fields
$\bf K$. Such fields have non-archimedean norms and their
characteristics may be either zero such as for $\bf Q_p$ or for its
finite algebraic extension, or positive characteristics $char ({\bf
K})=p>0$ such as ${\bf F_p}(\theta )$ of Laurent series over a
finite field ${\bf F_p}$ with $p$ elements and an indeterminate
$\theta $, where $p>1$ is a prime number \cite{wei}. Multiplicative
norms in such fields $\bf K$ satisfy stronger inequality, than the
triangle inequality, $|x+y|\le \max (|x|, |y|)$ for each $x, y\in
\bf K$. Non-archimedean fields are totally disconnected and balls in
them are either non-intersecting or one of them is contained in
another. In works \cite{byvo,evans}-\cite{evans4,khrkoz} stochastic
processes on spaces of functions with domains of definition in a
non-archimedean linear space and with ranges in the field of real
$\bf R$ or complex numbers $\bf C$ were considered. Different
variants of non-archimedean stochastic processes are possible
depending on a domain of definition, a range of values of functions,
values of measures in either the real field or a non-archimedean
field \cite{luijmms05,ludkhr}, a time parameter may be real or
non-archimedean and so on. That is, depending on considered problems
different non-archimedean variants arise.
\par Stochastic processes with values in non-archimedean spaces
appear while their studies for non-archimedean Banach spaces,
totally disconnected topological groups and manifolds
\cite{lusmfn2006}-\cite{luijmms04}. Very great importance have also
branching processes in graphs \cite{aigner,gihsko,hentor}. For
finite or infinite graphs with finite degrees of vertices there is
possible to consider their embeddings into $p$-adic graphs, which
can be embedded into locally compact fields. That is, a
consideration of such processes reduces to processes with values in
either the field $\bf Q_p$ of $p$-adic numbers or ${\bf F_p}(\theta
)$.
\par In this article theorems about representations of
characteristic functionals of infinitely divisible distributions
with values in vector spaces over locally compact infinite fields
with non-trivial non-archimedean valuations are formulated and
proved. Special features of the non-archimedean case are elucidated.
Therefore, a part of definitions, formulations of theorems and their
proofs are changed in comparison with the classical case. All main
results of this paper are obtained for the first time.
\par There is also an interesting interpretation of stochastic
processes with values in $\bf Q_p^n$, for which a time parameter may
be either real or $p$-adic. A random trajectory in $\bf Q_p^n$ may
be continuous relative to the norm induced by the non-archimedean
valuation in $\bf Q_p$, but its trajectory in $\bf Q^n$ relative to
the usual metric induced by the real metric may be discontinuous.
This gives new approach to spasmodic or jump or discontinuous
stochastic processes with values in $\bf Q^n$, when the latter is
considered as embedded into $\bf R^n$. On the other hand, stochastic
processes with values in ${\bf F_p}(\theta )^n$ can naturally take
into account cyclic stochastic processes in definite problems.

\par {\bf 1. Notations and definitions.} Let $(\Omega ,{\cal A},P)$ -
be a probability space, where $\Omega $ is a space of elementary
events, ${\cal A}$ is a $\sigma $-algebra of events in $\Omega $,
$P: {\cal A}\to [0,1]$ is a probability. Denote by $\xi $ a random
vector (a random variable for $n=1$) with values in $\bf K^n$ such
that it has the probability distribution $P_{\xi }(A) = P( \{ \omega
\in \Omega : \xi (\omega )\in A \} )$ for each $A\in {\cal B}({\bf
K^n})$, where $\xi : \Omega \to \bf K^n$, $\xi $ is $({\cal A},{\cal
B}({\bf K^n}))$-measurable, that is, $\xi ^{-1}({\cal B}({\bf
K^n}))\subset \cal A$, where $\bf K$ is a locally compact infinite
field with a non-trivial non-archimedean valuation, $n\in \bf N$,
$\bf Q_p$ is the field of $p$-adic numbers, $1<p$ is a prime number.
Here $\bf K$ is either a finite algebraic extension of the field
$\bf Q_p$ or the field $\bf Q_p$ itself for $char ({\bf K})=0$, or
${\bf K}={\bf F_p}(\theta )$ for $char ({\bf K})=p>1$, ${\cal
B}({\bf K^n})$ is the $\sigma $-algebra of all Borel subsets in $\bf
K^n$. Random vectors $\xi $ and $\eta $ with values in $\bf K^n$ are
called independent, if $P( \{ \xi \in A, \eta \in B \} )=P( \{ \xi
\in A \} ) P( \{ \eta \in B \} )$ for each $A, B\in {\cal B}({\bf
K^n})$.
\par A random vector (a random variable) $\xi $ is called infinitely
divisible, if
\par $(1)$ for each $m\in \bf N$ there exist random vectors (random
variables) $\xi _1,...,\xi _m$ such that $\xi = \xi _1+...+\xi _m$
and the probability distributions of $\xi _1,...,\xi _m$ are the
same.
\par If $\xi =\xi (t)=\xi (t,\omega )$ is a stochastic process with the
real time, $t\in T$, $T\subset \bf R$, then it is called infinitely
divisible, if Condition $(1)$ is satisfied for each $t\in T$.
Introduce the notation $B(X,x,R):= \{ y\in X: \rho (x,y)\le R \} $
for the ball in a metric space $(X,\rho )$ with a metric $\rho $,
$0<R<\infty $, $\xi _j(t)$ are stochastic processes, $j=1,...,m$.

\par {\bf 2. Lemma.} {\it If $\xi $ and $\eta $ are two independent
random vectors with values in $\bf K^n$ with probability
distributions $P_{\xi }$ and $P_{\eta }$, then $\xi +\eta $ has the
probability distribution $P_{\xi +\eta } (A)=\int_{\bf K^n} P_{\xi
}(A-dy)P_{\eta }(dy)$ for each $A\in {\cal B}({\bf K^n})$.}
\par {\bf Proof.} Since $\xi $ and $\eta $ are independent,
then $P(\{ \omega \in \Omega : \xi (\omega )\in C, \eta (\omega )\in
B \} )= P(\{ \omega \in \Omega : \xi (\omega )\in C \} ) P( \{
\omega \in \Omega : \eta (\omega )\in B \} )$ for each $C, B\in
{\cal B}({\bf K^n})$. Therefore, $P(\{ \xi +\eta \in A \} )= P( \{
\xi \in A-y, \eta =y, y\in {\bf K^n} \} )$ for each $A\in {\cal
B}({\bf K^n})$, consequently, $P_{\xi +\eta } (A)=\int_{\bf K^n}
P_{\xi }(A-dy)P_{\eta }(dy)$.
\par This means that $P_{\xi +\eta }=P_{\xi }* P_{\eta }$
is the convolution of measures $P_{\xi }$ and $P_{\eta }$.
\par {\bf 3. Corollary.} {\it If $\xi $ is an infinitely divisible
random vector, then $P_{\xi }=P^{*m}_{\xi _1}$ for each $m\in \bf
N$, where $P^{*m}_{\eta }$ denotes the $m$-fold convolution $P_{\eta
}$ with itself.}
\par {\bf Proof.} In view of Lemma 2 and Definition 1
$P_{\xi }=P_{\xi _1}*P_{\xi _2+...+\xi _m}=...=P_{\xi _1}*P_{\xi
_2}*...*P_{\xi _m}$. On the other hand, $\xi _1,...,\xi _m$ have the
same probability distributions, hence $P_{\xi _1}*P_{\xi
_2}*...*P_{\xi _m}=P^{*m}_{\xi _1}$.
\par {\bf 4. Notes and definitions.}
Corollary 3 means, that the equality $P_{\xi }=P^{*m}_{\xi _1}$
implies the relation: $P_{\xi }(A)=\int_{\bf K^n}...\int_{\bf K^n}
P_{\xi _1}(A-dy_2)P_{\xi _2}(dy_2-dy_3)...P_{\xi
_{m-1}}(dy_{m-1}-dy_m) P_{\xi _m}(dy_m)$, where $A\in {\cal B}({\bf
K^n})$. In the case of $char ({\bf K})=p>1$ Corollary 3 means, that
for $m=kp$, where $k\in \bf N$, if $P_{\xi }( \{ 0 \} )=0$, then
$P_{\xi _1}(\{ y \} )=0$ for each singleton $y\in \bf K^n$, since
$P(\xi =0)\ge P( \xi _1=\xi _2=...=\xi _m)\ge P_{\xi _1}( \{ y \}
)^m$. It is the restriction on the atomic property of $P_{\xi }$ and
$P_{\xi _1}$.

\par For $p$-adic numbers $x=\sum_{k=N}^{\infty }x_kp^k$, where
$x_k\in \{ 0, 1,...,p-1 \} $, $N\in \bf Z$, $N=N(x)$, $x_N\ne 0$,
$x_j=0$ for each $j<N$, put as usually $ord_{\bf Q_p}(x)=N$ for the
order of $x$, thus its norm is $|x|_{\bf Q_p} =p^{-N}$. Define the
function $[x]_{\bf Q_p} :=\sum_{k=N}^{-1}x_kp^k$ for $N<0$,
$[x]_{\bf Q_p}=0$ for $N\ge 0$ on ${\bf Q_p}$. Therefore, the
function $[x]_{\bf Q_p}$ on $\bf Q_p$ is considered with values in
the segment $[0,1]\subset \bf R$.
\par For the field ${\bf F_p}(\theta )$ put
$|x|_{{\bf F_p}(\theta )}= p^{-N}$, where $N=ord_{{\bf F_p}(\theta
)}(x) \in \bf Z$, $x=\sum_{j=N}^{\infty }x_j\theta ^j$, $x_j\in \bf
F_p$ for each $j$, $x_N\ne 0$, $x_j=0$ for each $j<N$. Then we
define the mapping $[x]_{{\bf F_p}(\theta )} =x_{-1}/p$, where we
consider elements of ${\bf F_p}=\{ 0, 1,...,p-1 \} $ embedded into
$\bf R$, hence $[x]_{{\bf F_p}(\theta )}$ takes values in $\bf R$,
where $1/p \in \bf R$, $x_{-1}=0$ when $N=N(x)\ge 0$.
\par Consider a local field $\bf K$ as the vector space over the field
$\bf Q_p$, then it is isomorphic with $\bf Q_p^b$ for some $b\in \bf
N$, since $\bf K$ is a finite algebraic extension of the field $\bf
Q_p$. In the case of ${\bf K}={\bf F_p}(\theta )$ we take $b=1$. Put
\par $(i)$ ${\bf F} := {\bf Q_p}$ for $char ({\bf K})=0$ with
${\bf K}\supset \bf Q_p$, while \par $(ii)$ ${\bf F} := {\bf
F_p}(\theta )$ for $char ({\bf K})=p>1$ with ${\bf K}={\bf
F_p}(\theta )$.
\par Let $(x,y):=(x,y)_{\bf F} := \sum_{j=1}^bx_jy_j$ for $x, y\in \bf
F$, $x=(x_1,...,x_b)$, $x_j\in \bf F$; $(x,y)_{\bf K}
:=\sum_{j=1}^nx_jy_j$ for $x, y\in \bf K^n$, $x=(x_1,...,x_n)$,
$x_j\in \bf K$. \par Define the mapping $<q>_{\bf F} := 2\pi
[(e,q)]_{\bf F}$ for each $q\in \bf K$, which is considered in
$(e,q)$ as the element from $\bf F^b$, $<q>_{\bf F}: {\bf K}\to \bf
R$, where $e:= (1,...,1)\in \bf F^b$, particularly $e=1$ for $b=1$,
that is, either in $(i)$ ${\bf K}=\bf Q_p$ or in the case $(ii)$ for
${\bf K}={\bf F_p}(\theta )$. For the additive group $\bf K^n$ then
there exists the character $\chi _s(z) := \exp (i <(s,z)_{\bf
K}>_{\bf F})$ with values in the field of complex numbers $\bf C$
for each value of the parameter $s\in \bf K^n$, since
$s_j(z_j+v_j)=s_jz_j+s_jv_j$ for each $s_j, z_j, v_j\in \bf K$ and
$(s,z+v)_{\bf K}=(s,z)_{\bf K}+(s,v)_{\bf K}$, $[x+y]_{\bf F} -
[x]_{\bf F} - [y]_{\bf F} \in B({\bf F},0,1)$ for every $x, y\in \bf
F$, while $[x]_{\bf F}=0$ for each $x\in B({\bf F},0,1)$, where
$i=(-1)^{1/2}\in \bf C$. In particular, $\chi _0(z)=1$ for each
$z\in \bf K^n$ for $s=0$. The character is non-trivial for $s\ne 0$.
At the same time $\chi _s(z)= \prod_{j=1}^n\chi _{s_j}(z_j)$, where
$\chi _{s_j}(z_j)$ are characters of $\bf K$ as the additive group.
\par For a $\sigma $-additive measure $\mu : {\cal B}({\bf K^n})\to \bf C$
of a bounded variation the characteristic functional $\hat \mu $ is
given by the formula: ${\hat \mu }(s) := \int_{\bf K^n} \chi
_s(z)\mu (dz)$, where $s\in \bf K^n$ is the corresponding continuous
$\bf K$-linear functional on $\bf K^n$ denoted by the same $s$.
\par In general the characteristic functional of the measure $\mu $
is defined in the space $C^0({\bf K^n},{\bf K})$ of continuous
functions $f: {\bf K^n}\to \bf K$
\par ${\hat \mu }(f) := \int_{\bf K^n} \chi _{1}(f(z))\mu
(dz)$, where $1\in \bf K$.
\par Let $\mu $ be a $\sigma $-additive finite non-negative
measure on ${\cal B}({\bf K^n})$, $\mu ({\bf K^n})<\infty $.
Consider the class ${\cal C}_1= {\cal C}_1({\bf K})$ of continuous
functions $A=A_{\mu }: {\bf K^n}\to \bf R$, satisfying Conditions
$(F1-F3)$:
\par $(F1)$ $A(y+z)=A(y)+A(z)+2\pi \int_{\bf K^n}f_1(y,z;x)\mu (dx)$
for each $y, z\in \bf K^n$,
\par $(F2)$ $A(\beta y)=[\beta ]_{\bf F}A(y)+2\pi \int_{\bf K^n}
f_2(\beta ,(e, (y,x)_{\bf K})_{\bf F}) \mu (dx)$ for each $y\in \bf
K^n$, $\beta \in \bf F$, where either \par $(F3)$ if ${\bf F}=\bf
Q_p$ for $char ({\bf K})=0$, then $f_1: ({\bf K^n})^3\to \bf Z$ and
$f_2: {\bf Q_p}^2\to \bf R$ are locally constant continuous bounded
functions, $f_1(y,z;x)\in \bf Z$ and $f_2(\alpha ,\beta )p^{ -
N(\alpha ,\beta )}\in \bf Z$ for $N(\alpha ,\beta )<0$ take only
integer values, $N(\alpha ,\beta ):=\min (ord_{\bf Q_p}(\alpha ),
ord_{\bf Q_p}(\beta ))$; or
\par $(F4)$ if ${\bf F}={\bf F_p}(\theta )$ for $char ({\bf K})=p >0$,
then $f_1: ({\bf K^n})^3\to \bf R$ and $f_2: {\bf F}^2\to \bf R$ are
locally constant continuous bounded functions, $pf_1(y,z;x)\in \bf
Z$ and $p^2f_2(\alpha ,\beta )\in \bf Z$ for $N(\alpha ,\beta )<0$
take only integer values, $N(\alpha ,\beta ):=\min (ord_{{\bf
F_p}(\theta )}(\alpha ), ord_{{\bf F_p}(\theta )}(\beta ))$. While
$f_1(y,z;x)=0$ for $\max( |yx|_{\bf K}, |zx|_{\bf K})\le 1$, and
$f_2(\alpha ,\beta )=0$ for $\max( |\alpha |_{\bf F}, |\beta |_{\bf
F})\le 1$ in $(F3,F4)$.

\par Denote by ${\cal C}_2= {\cal C}_2({\bf K})$ the class of
continuous functions
$B=B_{\mu }: ({\bf K^n})^2\to \bf R$, satisfying Conditions
$(B1-B3)$:
\par $(B1)$ $B(y,z)=B(z,y)$ for each $y, z\in \bf K^n$,
where $B(y,y)$ is non-negative,
\par $(B2)$ $B(q+y,z)=B(q,z)+B(y,z)+ 2\pi \int_{\bf K^n}
f_1(q,y;x)<(z,x)_{\bf K}>_{\bf F} \mu (dx)$ for each $q, y, z\in \bf
K^n$,
\par $(B3)$ $B(\beta y,z)=[\beta ]_{\bf F}B(y,z)+2\pi \int_{\bf K^n}
f_2(\beta ,(e,(y,x)_{\bf K})_{\bf F}) <(z,x)_{\bf K}>_{\bf F} \mu
(dx)$, where $f_1$ and $f_2$ satisfy Condition either $(F3)$ or
$(F4)$ depending on the characteristic $char ({\bf K})$. \par For
$y=z$ we shall also write for short $B(y):= B(y,y)$.

\par {\bf 4.1. Lemma.} {\it If $\chi _s(x): {\bf F^n}\to \bf C$
is a character of the additive group of $\bf F^n$ as in Section 4,
$\mu : {\cal B}({\bf F^n})\to [0,\infty ]$ is the Haar measure such
that $\mu (B({\bf F^n},0,1))=1$. Then $\int_{B({\bf F ^n},0,p^k)}
\chi _s(x)\mu (dx)= J(s,k)$, where $J(s,k)= p^{kn}$ for $|s|\le
p^{-k}$, while $J(s,k)=0$ for $|s|\ge p^{1-k}$.}
\par {\bf Proof.} The Haar measure $\mu $ on ${\cal B}({\bf F^n})$ is
the product of the Haar measures $\mu _1$ on ${\cal B}({\bf F})$,
$\mu (dx)= \otimes_{j=1}^n \mu _j(dx_j)$, $\mu _j=\mu _1$.
Therefore, $\int_{B({\bf F^n},0,p^k)} \chi _s(x)\mu (dx)= \prod
_{j=1}^n\chi _{s_j}(x_j)\mu _j(dx_j)$, where $\chi _j=\chi _1$,
$\chi _{s_j}(x_j)$ is the character of $\bf F$. \par Consider $n=1$.
Then $K := \int_{B({\bf F},0,p^k)} \chi _s(x)\mu (dx)=\int_{B({\bf
F},y,p^k)} \chi _s(x-y)\mu (dx)$  for each $y\in B({\bf F},0,p^k)$.
Thus $K=\chi _s(-y)\int_{B({\bf F},0,p^k)} \chi _s(x)\mu (dx)$,
since $B({\bf F},0,p^k)=B({\bf F},y,p^k)$ for each $y\in B({\bf
F},0,p^k)$, while $\mu (A-y)= \mu (A)$ for each $A\in {\cal B}({\bf
F})$. Take $|s|_{\bf F}\ge p^{-k+1}$ and $|y|_{\bf F}=p^k$ such that
$[sy]_{\bf F}\ne 0$ is nonzero. Hence $K (1-\chi _s(-y))=0$, but
$\chi _s(-y) \ne 1$, consequently, $K=0$. \par On the other hand, if
$|sx|_{\bf F}\le 1$, then $\chi _s(x)=1$ and inevitably
$\int_{B({\bf F},0,p^k)} \chi _s(x)\mu (dx)=p^k$, when $|s|_{\bf
F}\le p^{-k}$ (see for comparison the case ${\bf F}=\bf Q_p$ in
Example 6 on page 62 \cite{vla3}).

\par {\bf 5. Theorem.} {\it Let $\{ \psi (v,y): v\in V \} $ be a family of
characteristic functionals of $\sigma $-additive non-negative
bounded measures on ${\cal B}({\bf K^n})$, where $V$ is a
monotonically decreasing sequence of positive numbers converging to
zero. Suppose that there exists a limit $g(y)=\lim_{v\downarrow 0}
(\psi (v,y)-1)/v$ uniformly in each ball $B({\bf K^n},0,R)$ for each
given $0<R<\infty $. Then in $\{ {\bf K^n}, {\cal B}({\bf K^n}) \} $
there exists a $\sigma $-additive non-negative bounded measure $\nu
$, functions $A(y)$ and $B(y)$, belonging to classes ${\cal C}_1$
and ${\cal C}_2$ respectively such that
\par $(i)$ $g(y)=iA(y)- B(y)/2 + \int_{\bf K^n}(\exp (i <(y,x)_{\bf
K}>_{\bf F}) - 1 - i<(y,x)_{\bf K}>_{\bf F} (1+|x|^2)^{-1}+
<(y,x)_{\bf K}>_{\bf F}^2(1+|x|^2)^{-1}/2) [(1+|x|^2)/|x|^2]\nu
(dx)$, $\nu \ge 0$, $\nu ( \{ 0 \} )=0$.}
\par {\bf Proof.} Let $\mu _v$ be a measure corresponding to the
characteristic functional $\psi (v,y)$. Put $\lambda _v(A):= v^{-1}
\int_A |z|^2/[1+|z|^2]\mu _v(dz)$ for each $A\in {\cal B}({\bf
K^n})$, where $|z| := \max_{1\le j\le n} |z_j|_{\bf K}$,
$z=(z_1,...,z_n)\in \bf K^n$, $z_j\in \bf K$ for every $j=1,...,n$.
We prove a weak compactness of the family of measures $\{ \lambda
_v: v\in V \} $. That is, we need to prove that $(i)$ there exists
$L=const >0$ such that $\sup_{v\in V} \lambda _v({\bf K^n})\le L$;
$(ii)$ $\lim_{R\to \infty } {\overline {\lim }}_{v\downarrow 0}
\lambda _v({\bf K^n} \setminus B({\bf K^n},0,R))=0$.
\par The topologically dual space ${\bf K^n}'$ of all
continuous $\bf K$-linear functionals on $\bf K^n$ is $\bf
K$-linearly and topologically isomorphic with ${\bf K^n}$, since
$n\in \bf N$. Since $\bf K$ is the locally compact field, then it is
spherically complete (see Theorems 3.15, 5.36 and 5.39 \cite{roo}).
Since $\bf K^n$ as the linear space over $\bf F$ is isomorphic with
$\bf F^{bn}$, then it is sufficient to verify a weak compactness
over the field $\bf F$, where either ${\bf F}=\bf Q_p$ for $char
({\bf K})=0$ with ${\bf K}\supset \bf Q_p$ and $b\in \bf N$, or
${\bf F}={\bf F_p}(\theta )$ for $char ({\bf K})=p>0$ with ${\bf
K}={\bf F_p}(\theta )$ and $b=1$. Indeed, apply the non-archimedean
variant of the Minlos-Sazonov theorem, due to which there exists the
bijective correspondence between characteristic functionals and
measures \cite{lusmfn2006}, where characteristic functionals are
weakly continuous (see also \S IV.1.2 and Theorem IV.2.2 about the
Minlos-Sazonov theorem on Hausdorff completely regular (Tychonoff)
spaces \cite{vatach}). They are positive definite on $({\bf K^n})'$
or $C^0({\bf K^n},{\bf K})$, when $\mu $ is non-negative; ${\hat \mu
}(0)=1$ for $\mu ({\bf K^n})=1$. In the considered case $\bf K^n$ is
a finite dimensional Banach space over $\bf K$. Since the
multiplication in $\bf K$ is continuous, then over $\bf Q_p$ this
gives the continuous mapping $f_0: ({\bf Q_p^b})^2\to \bf Q_p^b$.
The composition of $f_0$ with all possible $\bf K$-linear continuous
functionals $s: {\bf K^n}\to \bf K$ separates points in $\bf K^n$.
\par Let $|x|\le R_1$, where $0<R_1<\infty $ is an arbitrarily
given number. Due to conditions of this theorem for each $\delta
>0$ there exists $v_0=v_0(R_1,\delta )>0$ such that for each
$\epsilon >0$ there is satisfied the inequality:
\par $(1)$ $ - Re g(y)+\delta \ge \int_{B({\bf F^{bn}},0,\epsilon )}
[1-\cos <(y,x)_{\bf F} >_{\bf F}]|x|^{-2}\lambda _v(dx)$ for each
$0<v\le v_0$, since $e^{i\alpha }=\cos (\alpha ) + i \sin (\alpha
)$, $-Re (e^{i\alpha }-1)=1-\cos (\alpha )$ for each $\alpha \in \bf
R$, while $1+|x|^2\ge 1$ and $[1+|x|^2]|x|^{-2}\ge |x|^{-2}$.

\par If $\epsilon >1$ and $x\in {\bf F^{bn}}\setminus B({\bf F^{bn}},
0,{\epsilon })$, then from $|x|_{\bf F}>\epsilon $ it follows
$[1+|x|^2]|x|^{-2} =1+|x|^{-2}\ge 1$ and then for each $\delta
>0$ there exists $v_0>0$ such that for each $\epsilon >1$ and each
$0<v\le v_0$ there is satisfied the inequality:
\par $(2)$ $ - Re g(y) +\delta \ge
\int_{{\bf F^{bn}}\setminus B({\bf F^{bn}}, 0,{\epsilon })} (1- \cos
<(y,x)>_{\bf F}) \lambda _v(dx)$.
\par Integrate these inequalities by $y\in B({\bf Q_p^{bn}},
0,r)$ and divide on the volume (measure) $\mu (B({\bf F^{bn}},
0,r))$, where $\mu $ is the nonnegative Haar measure on $\bf F^{bn}$
such that $\mu (B({\bf F^{bn}}, 0,1))=1$, $\mu (B({\bf
F^{bn}},0,r))=r^{bn}$ for each $r=p^k$ with $k\in \bf Z$
\cite{bourbmh,wei}. Then from $(1)$ it follows: \par $(3)$ $-
r^{-bn} \int_{B({\bf F^{bn}}, 0,r)} Re g(y)\mu (dy) + \delta \ge
\int_{ B({\bf F^{bn}}, 0,r)} (\int_{B({\bf F^{bn}}, 0,\epsilon )}
|x|_{\bf F}^{-2} (1- \cos <(y,x)>_{\bf F})\lambda _v(dx)) \mu
(dy)r^{-bn}$. From $(2)$ we get:
\par $(4)$ $- r^{-bn} \int_{B({\bf F^{bn}}, 0,r)} Re
g(y)\mu (dy) + \delta \ge \int_{ B({\bf F^{bn}}, 0,r)} (\int_{{\bf
F^{bn}}\setminus B({\bf F^{bn}}, 0,\epsilon )} (1- \cos <(y,x)>_{\bf
F})\lambda _v(dx)) \mu (dy)r^{-bn}$. On the other hand, $\cos
(<(y,x)>_{\bf F}) = \cos (\sum_{j=1}^{bn} <x_j y_j>_{\bf F})$, since
$(y,x)= \sum_{j=1}^{bn} y_jx_j$, also $<a+b>_{\bf F}=<a>_{\bf F} +
<b>_{\bf F} + 2w \pi $ for each $a, b\in \bf F$, where $w$ is an
integer number, $w=w(a,b)\in \bf Z$. For the characters integrals
are known due to Lemma 4.1: $\int_{B({\bf F^{bn}},0,p^k)}\chi
_s(x)\mu (dx)= \prod_{j=1}^{bn}\int_{B({\bf F},0,p^k)}\chi
_{s_j}(x_j)\mu _j(dx_j)=J(s,k)$, where $J(s,k)=p^{kbn}$ for
$|s|_{\bf F}\le p^{-k}$, $J(s,k)=0$ for $|s|_{\bf F}\ge p^{-k+1}$.
Since $(y,x)=(x,y)$ and $\cos (\alpha )= Re (e^{i\alpha })$ for each
$\alpha \in \bf R$, then $\int_{B({\bf F^{bn}},0,p^k)}\cos
<(y,x)>_{\bf F}\mu (dy)=J(x,k)$, since $J(x,k)\in \bf R$. Take in
$(3,4)$ $r=p^k$, then
\par $(5)$ $- p^{-kbn} \int_{B({\bf F^{bn}}, 0,p^k)} Re
g(y)\mu (dy) + \delta \ge  (\int_{B({\bf F^{bn}}, 0,\epsilon )}
|x|_{\bf F}^{-2} (1- J(x,k)p^{-kbn} )\lambda _v(dx))$
\par $(6)$ $- p^{-kbn} \int_{B({\bf F^{bn}}, 0,p^k)} Re
g(y)\mu (dy) + \delta \ge \int_{{\bf F^{bn}}\setminus B({\bf
F^{bn}}, 0,\epsilon )} (1- J(x,k)p^{-kbn})\lambda _v(dx))$. Since
$J(x,k)p^{-kbn}=1$ for $|x|_{\bf F}\le p^{-k}$, while
$J(x,k)p^{-kbn}=0$ for $|x|_{\bf F}\ge p^{-k+1}$, then for $\epsilon
>p^{-k+1}$ with $k\in \bf Z$, where $p\ge 2$, we get $(1- J(x,k)p^{-kbn})=1$
for $p^{-k+1}\le |x|_{\bf F}\le \epsilon $, then
\par $- p^{-kbn} \int_{B({\bf F^{bn}}, 0,p^k)} Re
g(y)\mu (dy) + \delta \ge  (\int_{B({\bf F^{bn}}, 0,\epsilon
)\setminus B({\bf F^{bn}}, 0,p^{-k})} |x|_{\bf F}^{-2}\lambda
_v(dx))$
\\  $\ge \epsilon ^{-2}[\lambda _v(B({\bf F^{bn}}, 0,\epsilon ))-
\lambda _v(B({\bf F^{bn}}, 0,p^{-k})]$, hence \par $(7)$ $[\lambda
_v(B({\bf F^{bn}}, 0,\epsilon ))- \lambda _v(B({\bf F^{bn}},
0,p^{-k})] \le \epsilon ^2 [\delta -p^{-kbn}\int_{B({\bf F^{bn}},
0,p^k)} Re g(y)\mu (dy)]$. In particular, for $\epsilon _k=
p^{-k+2}$ with $\epsilon _k\le \epsilon $ and $k\to \infty $
Inequality $(7)$ is satisfied. The summation of both parts of
Inequality $(7)$ by such $k$ gives:
\par $(8)$ $\lambda _v(B({\bf F^{bn}}, 0,\epsilon ))\le
L_1\delta - \sum_{k=k_0}^{\infty }p^{-kbn-2k+4}\int_{B({\bf F^{bn}},
0,p^k)} Re g(y)\mu (dy)$, where $L_1= p^4\sum_{k=k_0}^{\infty }
p^{-2k}=p^{4-2k_0}/(1-p^{-2})$, $k_0\in \bf Z$ is fixed. At the same
time from $(6)$ it follows:
\par $(9)$ $- p^{-kbn} \int_{B({\bf F^{bn}}, 0,p^k)} Re
g(y)\mu (dy) + \delta \ge \lambda _v({\bf F^{bn}}\setminus B({\bf
F^{bn}}, 0,\epsilon ))$ for $\epsilon >p^{-k+1}$. Therefore, due to
Inequalities $(8,9)$ there exists $L=const>0$ such that $\lambda
_v({\bf F^{bn}})= \lambda _v(B({\bf F^{bn}},0,\epsilon ))+\lambda
_v({\bf F^{bn}}\setminus B({\bf F^{bn}},0,\epsilon ))\le L$, for
each $v\in (0,v_0]$, where $L=const >0$. \par Due to conditions of
this theorem the function $g(y)$ is continuous and $g(0)=0$,
consequently, for each $\delta >0$ there exists sufficiently small
$0< R_1=p^{k_1} <\infty $ such that $R_1^{-bn}|\int_{B({\bf
F^{bn}},0,R_1)} Re g(y)\mu (dy)|<\delta $. In view of Inequality
$(9)$ for each $\epsilon
> \max (p^{-k_1+1},1)$ there is satisfied the inequality
$\lambda _v({\bf F^{bn}} \setminus B({\bf F^{bn}},0,\epsilon
))<2\delta $ for each $v\in (0,v_0]$, consequently, the family of
measures $\{ \lambda _v: v\in V \} $ is weakly compact.
\par Choose a sequence $h_n\downarrow 0$ such that
$\lambda _{v_n}$ is weakly convergent to some measure $\nu $ on
${\cal B}({\bf K^n})$. Due to conditions of this theorem and using
the decomposition of $\exp $  into the series, we get the
inequality:
\par $(10)$ $[\psi (v,y)-1]/v = \int_{\bf K^n} (\chi _y(x)-1)
[1+|x|_{\bf K}^2] |x|_{\bf K}^{-2} \lambda _v(dx)$
\par $= i A_v(y) - B_v(y)/2 + \int_{\bf K^n} f(y,x)\lambda _v(dx)$,
\\ where $A_v(y) =\int_{\bf K^n} <(y,x)_{\bf K}>_{\bf F}
|x|_{\bf K}^{-2}\lambda _v(dx)$, $B_v(y)= \int_{\bf K^n} <(y,x)_{\bf
K}>_{\bf F}^2|x|_{\bf K}^{-2}\lambda _v(dx)$,
\par $(10')$ $f(y,x)= (\exp (i<(y,x)_{\bf K}>_{\bf F}) -1 -
i<(y,x)_{\bf K}>_{\bf F} [1+|x|_{\bf K}^2]^{-1} + <(y,x)_{\bf
K}>_{\bf F}^2[1+|x|_{\bf K}^2]^{-1}/2) [1+|x|_{\bf K}^2] |x|_{\bf
K}^{-2}$.
\par The multiplier $[1+|x|_{\bf K}^2] |x|_{\bf K}^{-2}$
is continuous and bounded for $|x|\ge R$, where $0<R<\infty $,
$<(y,x)_{\bf K}>_{\bf F}=0$ for $|y|_{\bf K} |x|_{\bf K}\le 1$,
hence the function $f(y,x)$ is continuous, it is bounded, when $y$
varies in a bounded subset in $\bf K^n$, while $x\in \bf K^n$.
Therefore, there exists $\lim_{k\to \infty } \int_{\bf K^n} f(y,x)
\lambda _{v_k} (dx) =\int_{\bf K^n} f(y,x)\nu (dx)$. The functions
$<(y,x)_{\bf K}>_{\bf F} |x|_{\bf K}^{-2}$ and $<(y,x)_{\bf K}>_{\bf
F}<(z,x)_{\bf K}>_{\bf F}|x|_{\bf K}^{-2}$ are locally constant by
$x$ for each given value of the parameters $y$ and $z$. These
functions are zero, when $|y|_{\bf K} |x|_{\bf K}\le 1$, that is,
they are defined in the continuous manner to be zero at the zero
point $x=0$. Since there exists the limit in the left hand side of
Inequality $(10)$, then there exist $\lim_{k\to \infty } A_{v_k}
(y)=A(y)=\int_{\bf K^n} <(y,x)_{\bf K}>_{\bf F} |x|_{\bf K}^{-2}\nu
(dx)$ and $\lim_{k\to \infty } B_{v_k}(y)=B(y)= \int_{\bf K^n}
<(y,x)_{\bf K}>_{\bf F}^2|x|_{\bf K}^{-2}\nu (dx)$. At the same time
$B(y)\ge 0$ for each $y\in \bf K^n$.
\par Substitute the measure $\nu (U)$ on $\nu (U\setminus \{ 0 \} )$ and
denote it by the same symbol, where $U\in {\cal B}({\bf K^n})$. Due
to the fact that $f(y,0)=0$, $<(y,0)_{\bf K}>_{\bf F}=0$, then for
such substitution of the measure the values of integrals $\int_{\bf
K^n} f(y,x)\nu (dx)$, $A(y)=\int_{\bf K^n} <(y,x)_{\bf K}>_{\bf
F}|x|_{\bf K}^{-2} \nu (dx)$ and $B(y,z)=\int_{\bf K^n} <(y,x)_{\bf
K}>_{\bf F}<(y,z)_{\bf K}>_{\bf F}|x|_{\bf K}^{-2} \nu (dx)$ do not
change.
\par It is known that $[\alpha +\beta ]_{\bf F}=[\alpha ]_{\bf F}+
[\beta ]_{\bf F} + v(\alpha ,\beta )$, where $v(\alpha ,\beta )\in
\bf Z$ for ${\bf F}=\bf Q_p$, $p v(\alpha ,\beta )\in \bf Z$ for
${\bf F}={\bf F_p}(\theta )$, $0\le [\alpha ]_p\le 1$ for each
$\alpha , \beta \in \bf F$. Also $[\alpha \beta ]_{\bf F}=[\alpha
]_{\bf F}[\beta ]_{\bf F} + u(\alpha ,\beta )$, where $p^{-N(\alpha
,\beta )}u(\alpha ,\beta )\in \bf Z$ for ${\bf F}=\bf Q_p$,
$p^2u(\alpha ,\beta )\in \bf Z$ for ${\bf F}={\bf F_p}(\theta )$,
since $[\alpha ]_{\bf Q_p}[\beta ]_{\bf Q_p}=\sum_{k=N(\alpha
)}^{-1}\sum_{l=N(\beta )}^{-1}\alpha _k\beta _lp^{k+l}$, and
$[\alpha \beta ]_{\bf Q_p} =\sum_{N(\alpha )\le k, N(\beta )\le l,
k+l\le -1} \alpha _k\beta _lp^{k+l}$, where $\alpha
=\sum_{k=N(\alpha )}^{\infty }\alpha _kp^k\in \bf Q_p$, $\alpha
_k\in \{ 0, 1,...,p-1 \} $ for each $k\in \bf Z$, $\alpha _{N(\alpha
)}\ne 0$, while $[\alpha ]_{{\bf F_p}(\theta )}[\beta ]_{{\bf
F_p}(\theta )}=\alpha _{-1}\beta _{-1}p^{-2}$, and $[\alpha \beta
]_{{\bf F_p}(\theta )} =\sum_{N(\alpha )\le k, N(\beta )\le l, k+l=
-1} \alpha _k\beta _lp^{-1}$, where $\alpha =\sum_{k=N(\alpha
)}^{\infty }\alpha _k\theta ^k\in {\bf F_p}(\theta )$, $\alpha _k\in
{\bf F_p}$ for each $k\in \bf Z$, $\alpha _{N(\alpha )}\ne 0$
\cite{vla3,wei}. At the same time $[\alpha ]_{\bf F}=0$, when
$|\alpha |_{\bf F}\le 1$, hence $v(\alpha , \beta ) =0$ and
$u(\alpha ,\beta )=0$ for $\max (|\alpha |_{\bf F},|\beta |_{\bf
F})\le 1$. Then
\par $(11)$ $<(y+z,x)_{\bf K}>_{\bf F}=<(y,x)_{\bf K}>_{\bf F} +
<(z,x)_{\bf K}>_{\bf F} + 2\pi f_1(y,z;x)$, where $f_1\in \bf Z$ for
${\bf F}=\bf Q_p$, $pf_1\in \bf Z$ for ${\bf F}={\bf F_p}(\theta )$.
Since $<(y,x)_{\bf K}>_{\bf F}$ is locally constant and $0\le
[\alpha ]_{\bf F}\le 1$ for each $\alpha \in \bf F$, then there is
the inequality $ - 2\le f_1(y,z;x)\le 1$ for each $x, y, z\in \bf
K^n$ in $(11)$. On the other hand,
\par $(12)$ $<(\beta y,x)_{\bf K}>_{\bf F}=[\beta ]_{\bf F}
<(y,x)_{\bf K}>_{\bf F} + 2\pi f_2(\beta ,(e,(y,x)_{\bf K})_{\bf
F})$, where $f_2(\alpha , \beta )=u(\alpha ,\beta )$ for each
$\alpha ,\beta \in \bf F$, since $\bf F$ is naturally embedded into
$\bf K$ and $\beta (e,(y,x)_{\bf K})_{\bf F}=(e,(\beta y,x)_{\bf
K})_{\bf F}$. Since $[\alpha ]_{\bf F}\in [0,1]$ for each $\alpha
\in \bf F$, then $-1\le f_2(\alpha ,\gamma )\le 1$ for each $\alpha
\in \bf F$ and $\gamma =(e,(y,x)_{\bf K})_{\bf F}\in \bf F$ in
$(12)$. In view of the continuity and the locally constant behavior
of $<(y,x)_{\bf K}>_{\bf F}$ from this the continuity and local
constantness of $f_1$ and $f_2$ follow. Thus, $f_1$ and $f_2$
satisfy Conditions $(F3,F4)$ depending on $char ({\bf K})$.
Therefore, from $(11,12)$ we get the properties:
\par $(13)$ $A(y)=\int_{\bf K^n} <(y,x)_{\bf K}>_{\bf F}
|x|_{\bf K}^{-2}\nu (dx)$ and
\par $(14)$ $B(y,z)= \int_{\bf K^n} <(y,x)_{\bf K}>_{\bf F}<(z,x)_{\bf
K}>_{\bf F}|x|_{\bf K}^{-2}\nu (dx)$ with the measure $|x|_{\bf
K}^{-2}\nu (dx)$ here instead of the measure $\mu $ in $(F1-F4)$,
$(B1-B3)$. By the construction given above the measures in the
definitions of $A$ and $B$ are nonnegative and the functions in
integrals are nonnegative, then $A(y)$ and $B(y,z)$ take nonnegative
values. \par As the metric space $\bf K^n$ is complete separable and
hence is the Radon space (see Theorem 1.2 \cite{dal}), that is, the
class of compact subsets approximates from below each $\sigma
$-additive nonnegative finite measure on the Borel $\sigma $-algebra
${\cal B}({\bf K^n})$. In view of the finiteness and the $\sigma
$-additivity of the nonnegative measure $|x|_{\bf K}^{-2}\nu (dx)$
on ${\bf K^n}\setminus B({\bf K^n},0,1/|y|_{\bf K})$ for $|y|_{\bf
K}>0$, $<(y,x)_{\bf K}>_{\bf F}=0$ for $|(x,y)_{\bf K}|\le 1$ and
due to continuity and boundedness of the functions in integrals we
have that the mappings $A(y)$ and $B(y,z)$ are continuous.

\par {\bf 6. Corollary.} {\it Let the conditions of Theorem 5
be satisfied and there exists $J := \int_{\bf K^n} |x|_{\bf
K}^{-2}\nu (dx)<\infty $. Then $A(y)= - i(\partial \phi (\beta
,y)/\partial \beta )|_{\beta =0}$ and $B(y)= - (\partial ^2 \phi
(\beta ,y)/\partial \beta ^2)|_{\beta =0}$, where $\phi (\beta ,y)=
\int_{\bf K^n} \exp (i <(y,x)>_{\bf F}\beta )|x|_{\bf K}^{-2}\nu
(dx)$, $-1< \beta <1$.}
\par {\bf Proof.} In view of Theorem 5 there exist $A(y)$
and $B(y)$. At the same time the measure $\nu $ is nonnegative as
the weak limit of a weakly converging sequence of nonnegative
measures, consequently, the measure $\mu (dx):= |x|_{\bf K}^{-2}\nu
(dx)$ is nonnegative. In view of the supposition of this Lemma $0\le
\mu ({\bf K^n})=J<\infty $. If $J=0$, then $A(y)=0$, $B(y)=0$ and
$\phi (\beta ,y)=0$, then the statement of this Lemma is evident.
Therefore, there remains the case $J>0$. Consider the random
variable $\zeta := <(y,\eta )_{\bf K}>_{\bf F}$ with values in $\bf
R$, where $\eta $ is a random vector in $\bf K^n$ with the
probability distribution $P(dx) := J^{-1} |x|_{\bf K}^{-2}\nu (dx)$,
where $y\in \bf K^n$ is the given vector.
\par Then $\phi (\beta ,y)= J M\exp (i\beta \zeta )$, where $MX$ denotes
the mean value of the random variable $X$ with values in $\bf C$.
That is, $M\exp (i\beta \zeta )=\int_{\bf K^n}\exp (i\beta
<(y,x)_{\bf K}>_{\bf F}) P(dx)$. For $\zeta $ there exists the
second moment, since there exists $B(y)$ for each $y\in \bf K^n$. In
view of Theorem II.12.1 \cite{shir} about relations between moments
of the random variable and values of derivatives of their
characteristic functions at zero, we get the statement of this
Corollary.

\par {\bf 7. Theorem.} {\it Let the conditions of Theorem 5 be satisfied
and in addition measures $\mu _v(dx)$ posses finite moments of
$|x|_{\bf K}$ of the second order: $\int_{\bf K^n}|x|_{\bf K}^2\mu
_v(dx)<\infty $ $\forall v\in V$, then for $g(y)$ there is the
representation:
\par $(i)$ $g(y)=i {\tilde A}(y) - {\tilde B}(y)/2 + \int_{B({\bf
K^n},0,\epsilon )} (\exp (i <(y,x)_{\bf K}>_{\bf F}) - 1 - i
<(y,x)_{\bf K}>_{\bf F} + <(y,x)_{\bf K}>_{\bf F}^2/2) \eta (dx) +
\int_{{\bf K^n}\setminus B({\bf K^n},0,\epsilon )} (\exp (i
<(y,x)_{\bf K}>_{\bf F}) - 1) \eta (dx)$, where $\eta $ is a
nonnegative $\sigma $-additive measure on ${\cal B}({\bf K^n})$,
$\eta ( \{ 0 \} )=0$, ${\tilde A}(y)\in {\cal C}_1$, ${\tilde
B}(y,z)\in {\cal C}_2$.}
\par {\bf Proof.} Let $\eta _v(A) :=
v^{-1} \int_A |x|^2 \mu _v(dx)$, where $\{ \mu _v : v \} $ is the
family of measures corresponding to the characteristic functions
$\psi (v,y)$. At first we prove the weak compactness of the family
of measures $\{ \Psi _B(x)\eta _v(dx): v\in V \} $ for $B=B({\bf
K^n},0,R)$, $0<R<\infty $, where $\Psi _B(x)=1$ for $x\in B$, $\Psi
_B(x)=0$ for $x\notin B$, $\Psi _B(x)$ is the characteristic
function of the set $B$. Using the non-archimedean analog of the
Minlos-Sazonov theorem as in \S 5 we reduce the proof to the case of
measures on $\bf F^{bn}$. Take $0<R_1<\infty $. In view of the
conditions of this theorem for each $\delta
>0$ there exists $v_0=v_0 (R_1,\delta )>0$ such that for each
$\epsilon >0$ and each $0<v\le v_0$ there is accomplished the
inequality $- Re g(y) +\delta \ge \int_{\bf F^{bn}} [1-\cos
(<(y,x)>_{\bf F})] |x|^{-2} \eta _v(dx)$ due to the existence of
$\lim_{v\downarrow 0} [\psi (v,y)-1]/v=g(y)$ uniformly in the ball
of the radius $0<R_1<\infty $, $\forall y\in \bf F^{bn}:$ $|y|\le
R_1$. Integrate this inequality by $y\in B({\bf F^{bn}},0,r)$ and
divide on the volume $\mu (B({\bf F^{bn}},0,r))=r^{bn}$ for $r\in
\Gamma _{\bf F} := \{ |x|: x\ne 0, x\in {\bf F} \} = \{ p^k: k\in
{\bf Z} \} $, where $\mu $ is the Haar nonnegative nontrivial
measure on $\bf F^{bn}$. Then $- r^{-bn} \int_{B({\bf F^{bn}},0,r)}
Re g(y) \mu (dy) +\delta \ge r^{-bn}\int_{B({\bf
F^{bn}},0,r)}(\int_{\bf F^{bn}} [1-\cos <(y,x)>_{\bf F}]
|x|^{-2}\eta _v(dx)\mu (dy)\ge r^{-bn}\int_{B({\bf
F^{bn}},0,r)}(\int_{B({\bf F^{bn}},0,\epsilon )} [1-\cos
<(y,x)>_{\bf F}] |x|^{-2} \eta _v(dx)\mu (dy)$, since $\eta _v\ge 0$
and $\mu \ge 0$ are nonnegative measures. Since $\int_{B({\bf
F^{bn}},0,p^k)}\chi _s(x)\mu (dx)=J(s,k)$, where $J(s,k)=p^{kbn}$
for $|s|\le p^{-k}$, $J(s,k)=0$ for $|s|\ge p^{-k+1}$, then \par $-
p^{-bnk} \int_{B({\bf F^{bn}},0,p^k)} Re g(y) \mu (dy) +\delta \ge
\int_{B({\bf F^{bn}},0,\epsilon )} [1-p^{-bnk}J(x,k)] |x|^{-2} \eta
_v(dx)$. For $\epsilon
> p^{-k+1}$ we then get $[\eta _v(B({\bf F^{bn}},0,\epsilon )
-\eta _v (B({\bf F^{bn}},0,p^{-k}))]\le \epsilon ^2 [\delta -
p^{-bnk} \int_{B({\bf F^{bn}},0,p^k)}Re g(y)\mu (dy)]$. Then for
$\epsilon =p^{-k_0+2}$ and $\epsilon _k=p^{-k+2}\le \epsilon $,
$k\to \infty $ the summing of these inequalities leads to: $\eta
_v(B({\bf F^{bn}},0,\epsilon ))\le L_1\delta - \sum_{k=k_0}^{\infty
}p^{-kbn-2k+4}\int_{B({\bf F^{bn}}, 0,p^k)} Re g(y)\mu (dy)$, where
$L_1= p^{4-2k_0}/(1-p^{-2})$, $k_0\in \bf Z$ is fixed.
\par In view of the fact that the function
$g(y)$ is continuous and $g(0)=0$, then for each $\delta
>0$ there exists $0<R_1<\infty $ such that $R_1^{-bn}|\int_{B({\bf
F^{bn}},0,R_1)} Re g(y)\mu (dy)|<\delta $. Then for $\epsilon
=p^{-k_0+2}$ there is accomplished the inequality: $\eta _v(B({\bf
F^{bn}},0,\epsilon ))<2L_1 \delta $ for each $v\in (0,v_0]$. Since
$\int_{{\bf K^n}\setminus B} \Psi _B(x)\eta (dx)=0$, then the family
of measures $\{ \Psi _B\eta _v: v\in V \} $ is weakly compact for
each given $0<R<\infty $, $B=B({\bf K^n},0,R)$.
\par Let $0<\epsilon <\infty $, then $|\int_{B({\bf
K^n},0,\epsilon )} <(y,x)_{\bf K}>_{\bf F}\nu (dx)|<\infty $ and
$|\int_{{\bf K^n}\setminus B({\bf K^n},0,\epsilon )} <(y,x)_{\bf
K}>_{\bf F}|x|^{-2} \nu (dx)|<\infty $, then $J_e:= \int_{\bf K^n}
(\exp (i<(y,x)_{\bf K}>_{\bf F}) -1 - i<(y,x)_{\bf K}>_{\bf F}
[1+|x|^2]^{-1} + <(y,x)_{\bf K}>_{\bf F}^2[1+|x|^2]^{-1}/2)
[1+|x|^2] |x|^{-2} \nu (dx)=\int_{\bf K^n} (\exp (i<(y,x)_{\bf
K}>_{\bf F}) -1 - i<(y,x)_{\bf K}>_{\bf F} [1+|x|^2]^{-1} +
<(y,x)_{\bf K}>_{\bf F}^2[1+|x|^2]^{-1}/2) \eta (dx)$, where $\eta
(A) := \int_A [1+|x|^2] |x|^{-2}\nu (dx)$ for each $A\in {\cal
B}({\bf K^n})$. The measure $\eta \ge 0$ is nonnegative, since $\nu
\ge 0$ is nonnegative. From $\nu ( \{ 0 \} )=0$ it follows that
$\eta ( \{ 0 \} )=0$. The measure $\eta (A)$ is finite for each
$A\in {\cal B} ({\bf K^n}\setminus B({\bf K^n},0,\epsilon ))$, when
$0<\epsilon <\infty $, since $\nu ({\bf K^n})<\infty $ and
$|x|>\epsilon $ for $x\in {\bf K^n}\setminus B({\bf K^n},0,\epsilon
)$. Therefore, $J_e= (\int_{B({\bf K^n},0,\epsilon )} + \int_{{\bf
K^n}\setminus B({\bf K^n},0,\epsilon )} (\exp (i<(y,x)_{\bf K}>_{\bf
F}) -1)\eta (dx) + \int_{\bf K^n}( - i<(y,x)_{\bf K}>_{\bf F} +
<(y,x)_{\bf K}>_{\bf F}^2/2)|x|^{-2} \nu (dx)$. At the same time
$\int_{\bf K^n}( - i<(y,x)_{\bf K}>_{\bf F} + <(y,x)_{\bf K}>_{\bf
F}^2/2)|x|^{-2} \nu (dx)= \int_{B({\bf K^n},0,\epsilon )}( -
i<(y,x)_{\bf K}>_{\bf F} + <(y,x)_{\bf K}>_{\bf F}^2/2) \eta (dx) -
\int_{B({\bf K^n},0,\epsilon )}( - i<(y,x)_{\bf K}>_{\bf F} +
<(y,x)_{\bf K}>_{\bf F}^2/2)[(1+|x|^2) -1]|x|^{-2} \nu (dx) +
\int_{{\bf K^n}\setminus B({\bf K^n},0,\epsilon )} ( - i<(y,x)_{\bf
K}>_{\bf F} + <(y,x)_{\bf K}>_{\bf F}^2/2)|x|^{-2} \nu (dx)$, hence
\par $(1)$ $g(y)= i {\tilde A}(y) - {\tilde B}(y)/2 + \int_{B({\bf
K^n},0,\epsilon )} (\exp (i<(y,x)_{\bf K}>_{\bf F}) - 1 -
i<(y,x)_{\bf K}>_{\bf F} + <(y,x)_{\bf K}>_{\bf F}^2/2) \eta (dx) +
\int_{{\bf K^n}\setminus B({\bf K^n},0,\epsilon )} (\exp
(i<(y,x)_{\bf K}>_{\bf F}) - 1)\eta (dx)$, where ${\tilde A}(y) =
A(y)+ \int_{ B({\bf K^n},0,\epsilon )} <(y,x)_{\bf K}>_{\bf F}\nu
(dx) - \int_{{\bf K^n}\setminus B({\bf K^n},0,\epsilon )}
<(y,x)_{\bf K}>_{\bf F}|x|^{-2} \nu (dx)$, ${\tilde B}(y) = B(y)+
\int_{ B({\bf K^n},0,\epsilon )} <(y,x)_{\bf K}>_{\bf F}^2\nu (dx) -
\int_{{\bf K^n}\setminus B({\bf K^n},0,\epsilon )} <(y,x)_{\bf
K}>_{\bf F}^2|x|^{-2} \nu (dx)$. \par Using the expressions for
$A(y)$ and $B(y,z)$ from the proof of Theorem 5, we get
\par $(2)$ ${\tilde A}(y) = \int_{ B({\bf K^n},0,\epsilon )}
<(y,x)_{\bf K}>_{\bf F}\nu (dx) + \int_{B({\bf K^n},0,\epsilon )}
<(y,x)_{\bf K}>_{\bf F}|x|^{-2} \nu (dx)$,
\par $(3)$ ${\tilde B}(y,z) = \int_{ B({\bf K^n},0,\epsilon )}
<(y,x)_{\bf K}>_{\bf F} <(z,x)_{\bf K}>_{\bf F}\nu (dx) +
\int_{B({\bf K^n},0,\epsilon )} <(y,x)_{\bf K}>_{\bf F}<(z,x)_{\bf
K}>_{\bf F}|x|^{-2} \nu (dx)$. Due to identities 5$(11,12)$ with the
measure $[1+|x|^{-2}]\Psi _B\nu (dx)$ here as the measure  $\mu $ in
\S 4, with $B=B({\bf K^n},0,\epsilon )$, where $\Psi _B(x)$ is the
characteristic function of the set $B$, $\Psi _B(x)=1$ for $x\in B$,
$\Psi _B(x)=0$ for $x\in {\bf K^n}\setminus B$, we get, that $\tilde
A$ and $\tilde B$ satisfy Conditions $(F1-F4)$ and $(B1-B3)$
respectively. Since the measures in the definition of $\tilde A$ and
$\tilde B$ are nonnegative and the functions in integrals are
nonnegative, then ${\tilde A}(y)$ and ${\tilde B}(y,z)$ take
nonnegative values. \par As the metric space $\bf K^n$ is complete
and separable, hence it is the Radon space (see Theorem 1.2
\cite{dal}), that is the class of compact subsets approximates from
below each $\sigma $-additive nonnegative finite measure on the
Borel $\sigma $-algebra ${\cal B}({\bf K^n})$. In view of the
finiteness and $\sigma $-additivity of the nonnegative measure
$[1+|x|^{-2}]\Psi _B\nu (dx)$ and the boundedness of the continuous
functions in integrals the mappings ${\tilde A}(y)$ and ${\tilde
B}(y,z)$ are continuous.
\par {\bf 8. Theorem.} {\it A characteristic function $\psi (y)$
of an infinitely divisible distribution in $\bf K^n$ has the form
$\psi (y)=\exp (g(y))$, where $g(y)$ is given by Formula $5(i)$. If
in addition distributions $\mu _v(dx)$ from Theorem 5 posses finite
moments $|x|_{\bf K}$ of the second order: $\int_{\bf K^n}|x|_{\bf
K}^2\mu _v(dx)<\infty $, then $g(y)$ is given by Formula $7(i)$.}
\par {\bf Proof.} Let $h_k=1/k$, $k\in \bf N$,
hence $g(y)=\lim_{k\to \infty } (\psi _k(y)-1)/(1/k)= \lim_{k\to
\infty } k(\psi _k(y)-1)=\ln \psi (y)$, ãäå $\psi _k(y)=\psi
(1/k,y)$, $\psi (y)=[\psi _k(y)]^k$. If fix $arg \psi (0)=0$ and
take such a continuous branch $arg \psi (y)$, then $\psi (y)=\exp
(g(y))$, where $g(y)$ is given by Theorem 5 or 7.

\par {\bf 9. Definitions.} Let there is a random function $\xi (t)$
with values in ${\bf K^n}$, $t\in T$, where $(T,\rho )$ is a metric
space with a metric $\rho $. Then $\xi (t)$ is called stochastically
continuous at a point $t_0$, if for each $\epsilon
>0$ there exists $\lim_{\rho (t,t_0)\to 0} P(|\xi (t)-\xi
(t_0)|>\epsilon )=0$. If $\xi (t)$ is stochastically continuous at
each point of a subset $S$ in $T$, then it is called stochastically
continuous on $S$.
\par If $\lim_{R\to \infty } \sup_{t\in S} P(|\xi (t)|>R)=0$,
then a random function $\xi (t)$ is called stochastically bounded on
$S$.
\par Let $T=[0,a]$ or $T=[0,\infty )$, $a>0$. A random process
$\xi (t)$ with values in $\bf K^n$ is called a process with
independent increments, if $\forall n$, $0\le t_1<...<t_n$: random
vectors $\xi (0)$, $\xi (t_1)-\xi (0)$,...,$\xi (t_n)-\xi (t_{n-1})$
are mutually independent. At the same time the vector $\xi (0)$ is
called the initial state (value), and its distribution $P(\xi (0)\in
B)$, $B\in {\cal B}({\bf K^n})$, is called the initial distribution.
A process with independent increments is called homogeneous, if the
distribution $P(t,s,B):= P(\xi (t+s)-\xi (t)\in B)$, $B\in {\cal
B}({\bf K^n})$, of the vector $\xi (t+s)-\xi (t)$ is independent
from $t$, that is, $P(t,s,B)=P(s,B)$ for each $t<t+s\in T$.

\par {\bf 10. Theorem.} {\it Let $\psi (t,y)$ be a characteristic
function of the vector $\xi (t+s)-\xi (s)$, $t>0$, $s\ge 0$, where
$\xi (t)$ is the stochastically continuous random process with
independent increments with values in $\bf K^n$. Then $\psi
(t,y)=\exp (tg(y))$, where $g(y)$ is given by Formula $5(i)$. If in
addition $|\xi (t)|_{\bf K}$ has the second order finite moments,
then the function $g(y)$ is written by Formula $7(i)$.}
\par {\bf Proof.} Let $\xi (t)$ be a homogeneous stochastically
continuous process with independent increments with values in $\bf
K^n$, where $t\in T\subset \bf R$. Let $t>s$, then $|\psi (t,y)-
\psi (s,y)|= |M\exp (i<(y,\xi (t))_{\bf K}>_{\bf F}) - M \exp (i
<(y,\xi (s))_{\bf K}>_{\bf F})| = |M(\exp (i<(y,\xi (t)-\xi
(s))_{\bf K}>_{\bf F}) -1)\exp (i<(y,\xi (s))_{\bf K}>_{\bf F})|\le
M|\exp (i<(y,\xi (t)-\xi (s))_{\bf K}>_{\bf F})-1|$. Therefore, from
the stochastic continuity of $\xi (t)$ it follows continuity of
$\psi (t,y)$ by $t$. In view of being homogeneous and independency
of increments the equalities are accomplished $\psi (t_1+t_2,y)=
M\exp (i <(y,\xi (t_1+t_2) - \xi (t_1))>_{\bf F} + i <(y,\xi
(t_1)-\xi (0))>_{\bf F})) = M\exp (i <(y,\xi (t_1)-\xi (0))>_{\bf
F})) M \exp (i <(y,\xi (t_2) - \xi (0))>_{\bf F})=\psi (t_1,y)\psi
(t_2,y)$ for each $t_1, t_2\in T$. On the other hand, a unique
continuous solution of the equation $f(v+u)=f(v)f(u)$ for each $v,
u\in \bf R$ has the form $f(v)=\exp (a v)$, where $a\in \bf R$.
Thus, $\psi (t,y)=\exp (tg(y))$, where $g(y)=\lim_{t\downarrow 0}
(\psi (t,y)-1)/t$. Applying Theorems 5 and 7, we get the statement
of this theorem.
\par {\bf 11. Remark.} Consider auxiliary random process
$\eta := [\xi ]_p$ with values in $\bf R^n$, where
$[(q_1,...,q_n)]_p:=([q_1]_p,...,[q_n]_p)$ for $q=(q_1,...,q_n)\in
\bf K^n$. If $\xi (t)$ is a homogeneous process with independent
increments, then such is also $\eta $. Let $a(t) := M\eta (t)$ is a
mean value, while $R(t,s):= M[(\eta (t)-a(t))^* (\eta (s)-a(s))]$ is
the correlation matrix, where $\eta =(\eta _1,...,\eta _n)$ is the
row-vector, $A^*$ denotes the transposed matrix $A$. For the process
with independent increments and finite moments of the second order
then $R(t,s)=B(\min (t,s))$, where the matrix $B(t)$ is symmetric
and nonnegative definite. If $\xi (t)$ is the homogeneous process
with independent increments, $\eta $ has the finite second order
moments, then as it is known $a(t)=at$, $R(t,s)=B\min (t,s)$, where
$a$ is the vector, $B$ is the symmetric nonnegative definite matrix
\cite{gihsko}.
\par {\bf 12. Theorem.} {\it Let $P$ and $Q$ be two nonnegative
finite $\sigma $-additive measures on the Borel $\sigma $-algebra
${\cal B}({\bf K^n})$, where $\bf K$ is a locally compact infinite
field with a nontrivial non-archimedean valuation, $n\in \bf N$. If
their characteristic functions are equal ${\hat P}(y)={\hat Q}(y)$
for each $y\in \bf K^n$, then $P(A)=Q(A)$ for each $A\in {\cal
B}({\bf K^n})$.}
\par {\bf Proof.} The metric space $\bf K^n$
is complete and separable, consequently, it is the Radon space, then
$P$ and $Q$ are Radon measures (see Theorem 1.2 \cite{dal}). Then
for each $\delta >0$ there exists the ball $B({\bf K^n},z,R)$,
$0<R<\infty $, $z\in \bf K^n$, such that $P({\bf K^n}\setminus
B({\bf K^n},z,R))<\delta $ and $Q({\bf K^n}\setminus B({\bf
K^n},z,R))<\delta $. \par For each ball $B({\bf K^n},z,R_1)$, $z\in
\bf K^n$, $0<R_1<\infty $, due to the Stone-Weierstrass theorem for
each $\epsilon >0$ and each continuous bounded function $f: {\bf
K^n}\to \bf R$ there exist $b_1,...,b_k\in \bf C$ and
$s_1,...,s_k\in \bf K^n$ such that $\sup_{x\in B({\bf K^n},z,R_1)}
|b_1\chi _{s_1}(x)+...+b_k\chi _{s_k}(x) - f(x)|<\epsilon $, where
$\chi _s(x)$ is the character, $k\in \bf N$, since the family of all
finite $\bf C$-linear combinations of characters forms the algebra
which is the subalgebra of the algebra of all continuous functions
on $B({\bf K^n},z,R_1)$, the complex conjugation preserves this
subalgebra, this subalgebra contains all complex constants and
separates points in $B({\bf K^n},z,R_1)$ (see Theorem IV.10
\cite{reedsim}).
\par The characteristic function $\Psi _{B({\bf K^n},z,R)}$ of the set
$B({\bf K^n},z,R)$ is continuous on $\bf K^n$, since $\bf K^n$ is
totally disconnected and the ball $B({\bf K^n},z,R)$ is clopen in
$\bf K^n$ (simultaneously open and closed). Take $z\in \bf K^n$,
$0<\delta _k<1/k$, $0<\epsilon _k<1/k$, $R=R(\delta _k)\le R(\delta
_{k+1})$ for each $k$. For an arbitrary vector $z_1\in \bf K^n$ with
$|z-z_1|_{\bf K^n} < R(\delta _1)$ take the function $\Psi
^{\epsilon }(x)=b_1\chi _{s_1}(x)+...+b_v\chi _{s_v}(x)$ such that
\par $\sup_{x\in B({\bf K^n},z_1,R_1)}|\Psi ^{\epsilon _k}(x) - \Psi
_{B({\bf K^n},z_1,R_1)}(x)|<\epsilon _k$. Then \par $\int_{\bf
K^n}\Psi ^{\epsilon _k}(x)P(dx)=\int_{\bf K^n}\Psi ^{\epsilon
_k}(x)Q(dx)$ and \par $\int_{\bf K^n}\Psi _{B({\bf K^n},z_1,R_1)}(x)
P(dx) = P(B({\bf K^n},z_1,R_1))$, \par $\int_{\bf K^n}\Psi _{B({\bf
K^n},z_1,R_1)}(x) Q(dx)=Q(B({\bf K^n},z_1,R_1))$. On the other hand,
\par $|\int_{\bf K^n}\Psi _{B({\bf K^n},z_1,R_1)}(x)P(dx) - \int_{\bf
K^n}\Psi _{B({\bf K^n},z_1,R_1)}(x)Q(dx)|$ \par  $\le |\int_{\bf
K^n}\Psi ^{\epsilon _k}(x)P(dx) - \int_{\bf K^n}\Psi _{B({\bf
K^n},z_1,R_1)}(x)P(dx)|$ \par  $+ |\int_{\bf K^n}\Psi ^{\epsilon
_k}(x)Q(dx) - \int_{\bf K^n}\Psi _{B({\bf K^n},z_1,R_1)}(x)Q(dx)| +
|\int_{\bf K^n}\Psi ^{\epsilon _k}(x)P(dx)$ \par $- \int_{\bf
K^n}\Psi ^{\epsilon _k}(x)Q(dx)| \le \epsilon _k(P({\bf K^n}) +
Q({\bf K^n}))$. The right hand side of the latter inequality tends
to zero while $k\to \infty $, consequently, $P(B({\bf
K^n},z_1,R_1))=Q(B({\bf K^n},z_1,R_1))$ for each ball $B({\bf
K^n},z_1,R_1)$ in $\bf K^n$, where $0<R_1<\infty $, $z_1\in {\bf
K^n}$, since $\lim_{k\to \infty }\delta _k=0$ and $P({\bf
K^n}\setminus B({\bf K^n},z,R(\delta _k))<\delta _k$, $Q({\bf
K^n}\setminus B({\bf K^n},z,R(\delta _k))<\delta _k$. Since balls
form the base of the topology in $\bf K^n$, then $P(A)=Q(A)$ for
each $A\in {\cal B}({\bf K^n})$.
\par {\bf 13. Theorem.} {\it Random vectors $\eta _1,...,\eta _k$
in $\bf K^n$ are independent if and only if
\par $(1)$ $M\exp (i<(y_1,\eta _1)_{\bf K}+...+(y_k,\eta _k)_{\bf
K}>_{\bf F})= M\exp (i <(y_1,\eta _1)_{\bf K}>_{\bf F})...M\exp
(i<(y_k,\eta _k)_{\bf K}>_{\bf F})$ for each $y_1,...,y_k\in \bf
K^n$.}
\par {\bf Proof.} From the independence of $\eta _1,...,\eta _k$
it follows the independence of $<(y_1,\eta _1)_{\bf K}>_{\bf
F},...,<(y_k,\eta _k)_{\bf K}>_{\bf F}$, consequently, there is
satisfied the Equality $(1)$, since $\exp (i<(y_1,\eta _1)_{\bf
K}+...+(y_k,\eta _k)_{\bf K}>_{\bf F}) =\exp (i<(y_1,\eta _1)_{\bf
K}>_{\bf F})...\exp (i<(y_k,\eta _k)_{\bf K}>_{\bf F})$.
\par Vice versa let $(1)$ be satisfied. Denote by
$P_{\eta _1,...,\eta _k}$ the mutual probability distribution of
random vectors $\eta _1,...,\eta _k$, by $P_{\eta _j}$ denote the
probability distribution of $\eta _j$. Then $\int_{\bf K^n} \exp
(i<(y_1,x_1)_{\bf K}+...+(y_k,x_k)_{\bf K}>_{\bf F})P_{\eta
_1,...,\eta _k}(dx) = M\exp (i<(y_1,\eta _1)_{\bf K}+...+(y_k,\eta
_k)_{\bf K}>_{\bf F}) = M\exp (i <(y_1,\eta _1)_{\bf K}>_{\bf
F})...M\exp (i<(y_k,\eta _k)_{\bf K}>_{\bf F}) =
\prod_{j=1}^k\int_{\bf K^n} \exp (i<(y_j,x_j)_{\bf K}>_{\bf
F})P_{\eta _j}(dx_j)$, where $x=(x_1,...,x_k)$, $y_1,...,y_k,
x_1,...,x_k\in \bf K^n$. Therefore, by Theorem 12 $P_{\eta
_1,...,\eta _k}(A_1\times ... \times A_k) = P_{\eta _1}(A_1)...
P_{\eta _k}(A_k)$ for each $A_1,...,A_k\in {\cal B}({\bf K^n})$,
consequently, $\eta _1,...,\eta _k$ are independent.

\par {\bf 14. Definitions.} A sequence of random vectors
$\xi _m$ in $\bf K^n$ is called convergent by the distribution to a
random vector $\xi $, if for each continuous bounded function $f:
{\bf K^n}\to \bf R$ there exists $\lim_{m\to \infty } Mf(\xi
_m)=Mf(\xi )$.
\par Let a metric space $(X,\rho )$ be given with a metric
$\rho $ and a $\sigma $-algebra of Borel subsets ${\cal B}(X)$.
\par The family of probability measures $ {\cal P} := \{ P_{\beta }:
\beta \in \Lambda \} $ on $(X,{\cal B}(X))$, where $\Lambda $ is a
set, is called relatively compact, if an arbitrary sequence of
measures from $\cal P$ contains a subsequence weakly converging to
some probability measure.
\par A family of probability measures $ {\cal P} := \{ P_{\beta }:
\beta \in \Lambda \} $ on $(X,{\cal B}(X))$ is called dense, if for
each $\epsilon >0$ there exists a compact subset $C$ in $X$ such
that $\sup_{\beta \in \Lambda } P_{\beta }(E\setminus C)\le \epsilon
$.
\par A sequence $ \{ P_m: m\in {\bf N} \} $ of probability measures
$P_m$ is called weakly convergent to a measure $P$ when $m\to \infty
$, if for each continuous bounded function $f: X\to R$ there exists
$\lim_{m\to \infty } \int_Xf(x)P_m(dx)=\int_Xf(x)P(dx)$.

\par {\bf 15. Theorem.} {\it A random vector $\xi $ in $\bf K^n$
is a limit by a distribution of sums ${\tilde \xi }_m :=
\sum_{k=1}^m\xi _{m,k}$ of independent random vectors with the same
probability distribution $\xi _{m,k}$, $k=1,...,m$, if and only if
$\xi $ is infinitely divisible.}
\par {\bf Proof.} If $\xi $ is infinitely divisible, then
for each $m\ge 1$ there exists independent random vectors with the
same distribution $\xi _{m,1},...,\xi _{m,k}$ such that the
probability distributions of $\xi $ and of the sum $(\xi
_{m,1}+...+\xi _{m,k})$ are the same.
\par Let now ${\tilde \xi }_m$ be a sequence of arbitrary
vectors converging by the distribution to $\xi $ when $m\to \infty
$. Take $k\ge 1$ and group the summands writing ${\tilde \xi }_{mk}$
in the form: ${\tilde \xi }_{mk} = \zeta _{m,1}+...+\zeta _{m,k}$,
where $\zeta _{m,1}= \xi _{mk,1}+...+\xi _{mk,m}$, ...,$\zeta
_{mk,k} = \xi _{mk,m(k-1)+1}+...+\xi _{mk,mk}$. Since the sequence
${\tilde \xi }_{mk}$ converges by the distribution to $\xi $ while
$m\to \infty $, then the sequence of the probability distributions
$P_{{\tilde \xi }_{mk}}$ of random vectors ${\tilde \xi }_{mk}$ is
relatively compact, consequently, due to the Prohorov Theorem (see
\S VI.25 \cite{hentor} or III.2.1 \cite{shir}) it is dense.
\par On the other hand, if $|{\tilde \xi }_{mk}|>R$, then due to
non-archimedeanity of the norm in $\bf K^n$ there exists $j$ such
that $|\zeta _{m,j}|>R$, consequently, $P(\zeta _{m,1}\in {\bf
K^n}\setminus B({\bf K^n},0,R))\le P({\tilde \xi }_{mk}\in {\bf
K^n}\setminus B({\bf K^n},0,R))$, since $\zeta _{m,j}$ are
independent and have the same probability distribution. Therefore,
$\{ P_{\zeta _{m,1}} : m\in {\bf N} \} $ is the dense family of
probability distributions. Then there exists the sequence $ \{ m_j:
j\in {\bf N} \} $ and random vectors $\eta _1,...,\eta _k$ such that
$\zeta _{m_j,l}$ converges by the distribution to $\eta _l$ for each
$l=1,...,k$ for $j\to \infty $. In view of the definition of
convergence by the distribution this means in particular, that for
each $b_1,...,b_k\in \bf K^n$ there exists $\lim_{j\to \infty }
M\exp (i <(b_1,\zeta _{m_j,1})_{\bf K}+...+(b_k,\zeta _{m_j,k})_{\bf
K}>_{\bf F}) = M\exp (i <(b_1,\eta _1)_{\bf K}+...+(b_k,\eta
_k)_{\bf K}>_{\bf F})$. In view of independency of random vectors
$\zeta _{m_j,1},...,\zeta _{m_j,k}$ there is satisfied the equality
$M\exp (i <(b_1,\zeta _{m_j,1})_{\bf K} +...+(b_k,\zeta
_{m_j,k})_{\bf K}>_{\bf F}) =M\exp (i <(b_1,\zeta _{m_j,1})_{\bf
K}>)...M\exp (i(b_k,\zeta _{m_j,k})_{\bf K}>_{\bf F})$, since $\exp
(i<y>_{\bf F})$ is the character of the additive group of the field
$\bf K$. Therefore, $\lim_{j\to \infty }M\exp (i <(b_1,\zeta
_{m_j,1})_{\bf K}+...+(b_k,\zeta _{m_j,k})_{\bf K}>_{\bf F}) = M\exp
(i <(b_1,\eta _1)_{\bf K}>_{\bf F})...M\exp (i<(b_k,\eta _k)_{\bf
K}>_{\bf F})$, thus, $M\exp (i <(b_1,\eta _1)_{\bf K}+...+(b_k,\eta
_k)_{\bf K}>_{\bf F})= M\exp (i <(b_1,\eta _1)_{\bf K}>_{\bf
F})...M\exp (i<(b_k,\eta _k)_{\bf K}>_{\bf F})$ for each
$b_1,...,b_k\in \bf K^n$. Then from Theorem 13 it follows, that the
random vectors $\eta _1,...,\eta _k$ are independent.
\par Since ${\tilde \xi }_{m_jk} = \zeta _{m_j,1} +...+ \zeta
_{m_j,k}$ converges by the distribution to $\eta _1+...+\eta _k$ and
${\tilde \xi }_{m_jk}$ converges by the distribution to $\xi $, then
$\xi $ is equal to $\eta _1+...+\eta _k$ by the distribution, since
$Mf(\xi )= \lim_{j\to \infty } Mf({\tilde \xi }_{m_jk})=\lim_{j\to
\infty } Mf(\zeta _{m_j,1} +...+ \zeta _{m_j,k})=Mf(\eta _1+...+\eta
_k)$ for each continuous bounded function $f: {\bf K^n}\to \bf R$.

\par {\bf 16. Particular cases of Theorem 10.}
1. If $A(y)= q <(a,y)_{\bf K}>_{\bf F}$, $B=0$, $\nu =0$, where
$a\in \bf K^n$ is some vector, $q=const
>0$, then $\psi (t,y)=\exp (i t q<(a,y)_{\bf K}>_{\bf F})$. The  the
random function $\eta (t)=<(\xi (t),y)_{\bf K}>_{\bf F}$ has the
form $\eta (t)=\eta (0) + tq$, where $\xi $ is the initial random
vector with values in $\bf K^n$. That is, $\eta (t)$ corresponds to
the uniform motion of the point in $\bf R$ with the velocity $q$.
\par In the case, when $A(y)= q (v,<y>_{\bf F})$, $B=0$, $\nu =0$,
where $v\in \bf R^n$ is a given vector, $0\le v_j\le 1$ for each
$j=1,...,n$, $v=(v_1,...,v_n)$, $q=const
>0$, then $\psi (t,y)=\exp (i t q(v,<y>_{\bf F}))$. Therefore,
the random variable $\eta (t)=(<\xi (t)>_{\bf F},<y>_{\bf F})_{\bf
R}$ has the form $\eta (t)=\eta (0) + tq$.
\par 2. It is possible to consider in formulas for
$A(y)$ and $B(y,z)$ in \S \S 5 and 7 in particular atomic measures,
denoting $\tilde A$ by $A$ and $\tilde B$ by $B$ here for the
uniformity, then there are the expressions of the form $\sum_j q_j
<(x_j,y)_{\bf K}>_{\bf F}$ and $\sum_jq_j<(x_j,y)_{\bf K}>_{\bf F}
<(x_j,z)_{\bf K}>_{\bf F}$, where $q_j=\nu ( \{ x_j \} )>0$ or
$q_j=\nu ( \{ x_j \} )|x_j|^{-2}>0$ depending on the considered
case, $x_j\ne 0$. In particular, there may be
$x_j=e_j=(0,...,0,1,0,...)\in \bf K^n$ with the unity on the $j$-th
place. These expressions may be transformed using conditions
$(F1-F4)$ or $(B1-B3)$ (see Formulas 5$(i,10,13,14)$ or 7$(i,1-3)$).
Then there are possible cases $A(y)=q<(a,y)_{\bf K}>_{\bf F}$,
$A(y)=(v,<y>_{\bf F})_{\bf R}$, $B(y,z)=\sum_{j=1}^n<s_jy_jz_j>_{\bf
F}$, $B(y,z)= \sum_{j=1}^nq_j<y_j>_{\bf F}<z_j>_{\bf F}$, where
$<y>_{\bf F}=(<y_1>_{\bf F},...,<y_n>_{\bf F})$, $y=(y_1,...,y_n)\in
\bf K^n$, $y_k\in \bf K$ for each $k$, $v\in \bf R^n$, $(*,*)_{\bf
R}$ is the scalar product in $\bf R^n$, $s_j\in \bf K$, $a\in {\bf
K^n}$. The consideration of the transition matrix $Y$ from one basis
in $\bf K^n$ to another or the matrix $X$ of transition from one
basis in $\bf R^n$ into another leads to the more general
expressions for $B(y,z)$ such as $B(y,z) = (b<y>_{\bf F},<z>_{\bf
F})_{\bf R}$, $B(y,z) = <(hy,z)_{\bf K}>_{\bf F}$, where $b$ is the
symmetric nonnegative definite $n\times n$ matrix with elements in
the field of real numbers $\bf R$, $h$ is the symmetric $n\times n$
matrix with elements in the locally compact field $\bf K$.
\par 3. If $A(y)=q <(a,y)_{\bf K}>_{\bf F}$, $B(y,z)=
<(hy,z)_{\bf K}>_{\bf F}$, where $a\in \bf K^n$, $h$ is the
symmetric $n\times n$ matrix with elements in the field $\bf K$, if
the correlation term $\int_{\bf K^n}f(y,x)\nu (dx)=0$ from \S 5 or
$\int_{B({\bf K^n},0,\epsilon )} (\exp (i <(y,x)_{\bf K}>_{\bf F}) -
1 - i <(y,x)_{\bf K}>_{\bf F} + <(y,x)_{\bf K}>_{\bf F}^2/2) \eta
(dx) + \int_{{\bf K^n}\setminus B({\bf K^n},0,\epsilon )} (\exp (i
<(y,x)_{\bf K}>_{\bf F}) - 1) \eta (dx)=0$ from \S 7 is zero, then
$\psi (t,y) = \exp (i t q<(a,y)_{\bf K}>_{\bf F} - t <(hy,y)_{\bf
K}>_p/2)$. Then $\xi (t)$ is one of the non-archimedean variants of
the Gaussian process.
\par 4. In the case, when $A(y)=(v,<y>_{\bf F})_{\bf
R}$, $B(y,z) = (b<y>_{\bf F},<z>_{\bf F})_{\bf R}$ (see paragraph
2), while the correlation term is zero, then $\psi (t,y) = \exp (i t
(v,<y>_{\bf F})_{\bf R} - t (b<y>_{\bf F},<y>_{\bf F})_{\bf R}/2)$
and again $\xi (t)$ is one of the analogs of the Gaussian process.
Though Gaussian processes in the non-archimedean case do not exist.
That is, we can satisfy a part of properties of the Gaussian type in
the non-archimedean case, but not all (see also \cite{luijmms05}).
\par 5. When $A=0$, $B=0$ (taking into account $(F1-F4)$ and
$(B1-B3)$; see Formulas 5$(i,10,13,14)$ or 7$(i,1-3)$), where $\nu $
is the purely atomic measure, concentrated at the point $z_0$, $\nu
( \{ z_0 \} )=q>0$, then $\psi (t,y)=\exp (qt (\exp (i <(y,z_0)_{\bf
K}>_{\bf F}) - 1)$. Therefore, $\xi (t)$ is the non-archimedean
analog of the Poisson process.
\par 6. If ${\tilde A}(y)=q <(a,y)_{\bf K}>_{\bf F}+ \int_{B({\bf
K^n},0,\epsilon )} <(y,x)_{\bf K}>_{\bf F}\eta (dx)$, ${\tilde
B}(y)= - \int_{B({\bf K^n},0,\epsilon )} <(y,x)_{\bf K}>_{\bf F}^2
\eta (dx)/2$, $\eta (B({\bf K^n},0,\epsilon ))<\infty $, then $g(y)
= i <(a,y)_{\bf K}>_{\bf F} + w \int_{\bf K^n} (\exp (i<(y,x)_{\bf
K}>_{\bf F})-1)\lambda (dx)$, where $\lambda $ is the probability
measure on $({\bf K^n},{\cal B}({\bf K^n}))$, $0<w=\eta ({\bf
K^n})<\infty $, $\eta (dx)=w\lambda (dx)$ (see Formulas 7$(i,1-3)$
and $(F1-F4)$, $(B1-B3)$). Therefore, $\psi (t,y)= \exp (i t q
<(a,y)_{\bf K}>_{\bf F}) \sum_{k=0}^{\infty } \exp (-wt) ((wt)^k/k!)
[\int_{\bf K^n} \exp (i <(y,x)_{\bf K}>_{\bf F})\lambda (dx)]^k$.
This expression of the characteristic function of the random process
$\xi (t) = \rho (t) + \xi _1 +...+\xi _{\zeta (t)}$, where $\rho
(t)$ is the random process in $\bf K^n$ with the characteristic
function $\exp (i t q <(a,y)_{\bf K}>_{\bf F})$, where $\xi
_1,...,\xi _k,...$ are independent random vectors in $\bf K^n$ with
the same probability distribution $\lambda (dx)$, $\zeta (t)$ is the
Poisson process with a parameter $w$ independent from $\rho , \xi
_1,...,\xi _k,...$. Then there arises the non-archimedean analog
$\xi (t)$ of the generalized Poisson process.
\par If ${\tilde A}(y)=(v,<y>_{\bf F})_{\bf R} + \int_{B({\bf
K^n},0,\epsilon )} <(y,x)_{\bf K}>_{\bf F}\eta (dx)$, where ${\tilde
B}(y)$ is the same as at the beginning of the given paragraph, then
$\psi (t,y)= \exp (i t (v,<y>_{\bf F})_{\bf R}) \sum_{k=0}^{\infty }
\exp (-wt) ((wt)^k/k!) [\int_{\bf K^n} \exp (i <(y,x)_{\bf K}>_{\bf
F})\lambda (dx)]^k$, where $\rho (t)$ has the characteristic
function $\exp (i t (v,<y>_{\bf F})_{\bf R})$.
\par {\bf 17. Remark.} Let a branching random process
is realized with values in the ring $\bf Z_p$ of integer $p$-adic
numbers or in the ring $B({\bf F_p}(\theta ),0,1)$, denote it by
$\bf B$. In the particular case of the uniform distribution
$|x|^{-2}\nu (dx)$ in $\bf B$ the measure $\nu $ is proportional to
the Haar measure $\mu $, $|x|^{-2}\nu (dx) = q\mu (dx)$, where
$q>0$, $\mu ({\bf B})=1$, $\nu ({\bf F}\setminus {\bf B})=0$, ${\bf
K}={\bf F}=\bf Q_p$ or ${\bf K}={\bf F}={\bf F_p}(\theta )$
respectively here, $n=1$. Then it is possible to calculate $A(y)$
and $B(y)$. In view of \S 5 in this particular case $A(y)=q
\int_{\bf B} <yx>_{\bf F} \mu (dx)$ and $B(y)= q \int_{\bf B}
<yx>_{\bf F}^2 \mu (dx)$. If $y=0$, then $A(0)=0$ and $B(0)=0$,
therefore, consider the case $y\ne 0$. The function $<yx>_{\bf F}$
takes the zero value when $|yx|_{\bf F}\le 1$ and is different from
zero when $|x|_{\bf F}> 1/|y|_{\bf F}$.
\par In the considered case the support of the measure
$\nu $ is contained in $\bf B$, then $A(y)$ and $B(y)$ are equal to
zero when $|y|_{\bf F}\le 1$. But the Haar measure is invariant
relative to shifts $\mu (A+z)=\mu (A)$ for each Borel subset in $\bf
F$ with the finite measure $\mu (A)<\infty $ and each $z\in \bf F$.
Moreover, $\mu (zdx)=|z|_{\bf F}\mu (dx)$, where $|z|_{\bf
F}=p^{-ord_{\bf F}(z)}$ (see \cite{wei}). Then $A(y)=q \int_{z\in
{\bf F}, |y|_{\bf F}\ge |z|_{\bf F}>1} <z>_{\bf F} \mu (dz)/|y|_{\bf
F}$ and $B(y)=q \int_{z\in {\bf F}, |y|_{\bf F}\ge |z|_{\bf F}>1}
<z>_{\bf F}^2 \mu (dz)/|y|_{\bf F}$, where $|y|_{\bf F}>1$. At the
same time $z=\sum_{k=N(x)}^{\infty }z_kp^k$ for ${\bf F}=\bf Q_p$ or
$z=\sum_{k=N(x)}^{\infty }z_k\theta ^k$ for ${\bf F}={\bf
F_p}(\theta )$, where $N(z)= ord_p(z)$, $z_k\in \{ 0,1,...,p-1 \} $
or $z_k\in {\bf F_p}$. If $\nu (dx)=q\mu (dx)$, then $A(y)=q
|y|_{\bf F} \int_{z\in {\bf F}, |y|_{\bf F}\ge |z|_{\bf F}>1}
<z>_{\bf F} |z|_{\bf F}^{-2}\mu (dz)$ and $B(y)=q |y|_{\bf F}
\int_{z\in {\bf F}, |y|_{\bf F}\ge |z|_{\bf F}>1} <z>_{\bf F}^2
|z|_{\bf F}^{-2}\mu (dz)$. These integrals are expressible in the
form of finite sums, since $\mu (B({\bf F},x,p^k))=p^k$ for each
$k\in \bf Z$ and $z\in \bf F$, where the functions in the integrals
are locally constant.
\par The measure $\nu $ is Borelian, $\nu : {\cal B}({\bf K^n})\to
[0,\infty )$, therefore each its atom may be only a singleton. More
generally (see Formulas 5$(i,13,14)$), if $\nu =\nu _1+ \nu _2$,
where $\nu _2$ is the atomic measure, while $\nu _1(dx)=f(x)\mu
(dx)$, where $f(x)=g(|x|_{\bf F},<x>_{\bf F})$, $g: {\bf R^2}\to
[0,\infty )$ is a continuous function, then \par $A(y)=\sum_j
<yx_j>_{\bf F}|x_j|^{-2}\nu _2( \{ x_j \} ) + \int_{\bf F} <yx>_{\bf
F} f(x) |x|_{\bf F}^{-2} \mu (dx)$,
\par $B(y) = \sum_j <yx_j>_{\bf F}^2 |x_j|^{-2}\nu _2( \{ x_j \} ) +
\int_{\bf F} <yx>_{\bf F}^2 f(x)|x|_{\bf F}^{-2} \mu (dx)$, \\
where $ \{ x_j \} $ are atoms of the measure $\nu _2$, $\nu _2 ( \{
x_j \} )>0$, each $x_j\ne 0$ is nonzero. At the same time integrals
by the Haar measure $\mu $ on $\bf F$ with functions $<yx>_{\bf F}
f(x) |x|_{\bf F}^{-2}$ and $<yx>_{\bf F}^2 f(x) |x|_{\bf F}^{-2}$,
where $f(x)=g(|x|_{\bf F}, <x>_{\bf F})$, are expressible in the
form of series, since $|x|_{\bf F}$ and $<x>_{\bf F}$ are locally
constant, hence $f$ is locally constant.

\par Department of Applied Mathematics,
\par Moscow State Technical University MIREA,
\par Av. Vernadsky 78, Moscow 119454, Russia

\par e-mail: sludkowski@mail.ru

\end{document}